\documentclass[10pt]{article}

\usepackage{amssymb}
\usepackage{amsmath}
\usepackage{amsthm} 
\usepackage{mathrsfs}                                                           
\usepackage{slashbox}  
\usepackage[all]{xy}
\usepackage{wrapfig}

\def\2{{1\over 2}}

\newcommand{\rf}[1]{(\ref{#1})}

\newcommand{\ud}{\mathrm{d}}
\renewcommand{\t}{\tilde}
\newcommand{\p}{\partial}

\title{\bf{Unitary representations of a loop $ax+b$ group, \\
Wiener measure and $\Gamma$-function}}
\author{Anton M. Zeitlin\footnote{anton.zeitlin@yale.edu, http://math.yale.edu/$\sim$az84 http://www.ipme.ru/zam.html}\\   
Department of Mathematics,\\
Yale University,\\
442 Dunham Lab, 10 Hillhouse Avenue,\\
New Haven, CT 06511}

\begin{document}
\maketitle
\begin{abstract}
We construct a family of irreducible unitary representations of the loop affine group of a line ($ax+b$ group) 
with central extension on the Hilbert space of square integrable functions with respect to the Wiener measure. We relate the matrix coefficients of the elements of the loop $ax+b$ group to the loop analogue of the $\Gamma$-function. 
\end{abstract}
\section{Introduction}
There are two low-dimensional non-compact groups of high importance in representation theory: Heisenberg group and $ax+b$ group, the latter being an affine group of a real line. The representation theory of these groups reflects major difficulties arising in the study of noncompact real Lie groups of higher dimensions (e.g. irreducible representations are infinite-dimensional). The representation theory of the Heisenberg group (see e.g. \cite{taylor}) is fairly standard and  governed by Stone-von Neumann theorem, while the representations of $ax+b$ group (which we will refer to below as $G$) are more subtle and interesting due to their relation to the special function theory \cite{vilenkin}. 

The representations of the loop counterpart of the Heisenberg group has been studied for a long time and they become an important element in the representation theory of affine Lie algebras and the theory of vertex algebras.

In this article we study the representation theory of the loop counterpart of the $ax+b$ group ($\Omega G$) and its central extension ($\hat G$). We construct the unitary representations of  $\Omega G$, which naturally generalize the irreducible representations of $G$. The irreducible unitary representations of $G$ divide into three classes (up to equivalence). The first two classes contain infinite-dimensional representations and there is only one representation in each class, however the third class consists of 1-dimenional representations, labeled by the real parameter.

 The infinite dimensional representations of $G$ can be realized in the Hilbert space $L^2(\mathbb{R}_{+},\frac{dx}{x})$.  The construction of these representations is given by the canonical construction 
of induced representations (see Section 2 or \cite{gn},\cite{vilenkin} for more details). 
However, in the case of loop $ax+b$ group the 
situation seems to be more subtle. We take the $L^2$ space with respect to 
Wiener measure as a representation space. The unitary representations, which we construct, are labeled 
by a certain function. We prove that certain representations in this class (e.g. when this function is  
constant) are irreducible. We do not classify all unitary representations of $\hat G$ in this article. 
However, we make a conjecture, which we support by examples, that all the 
unitary representations of $\hat G$ divide into three classes, as well as in the finite dimensional case.  

It is worthy to note of another construction of irreducible unitary representations, started by Gelfand and coworkers 
in \cite{ggv}, where the authors study unitary representations of more general current groups without central 
extension. These representations are constructed via the so-called multiplicative integral of representations. It would be interesting to find the relation of our representations and the ones studied by the authors of \cite{ggv}.

An interesting issue, which we address in the second part of the article is based on the mentioned relations of 
the representations of group $G$ and the theory of special functions. It is known that the representations of $G$ are related to $\Gamma$- and $B$- functions \cite{vilenkin}. Therefore, we may hope to construct their loop counterparts while studying the associated loop group. In this article we make a first step in this direction. Namely, we use the fact that one can relate via the bilateral Laplace transform (which is a modification of the standard Fourier transform) the 
action of the group element of representation of $G$ in the representation space with the integral operator with the kernel, expressed as 
$\Gamma$-function. For the classical Wiener measure, one can define a certain analogue of the Fourier transform. Then,  
considering the action of the group elements in the representation space of $\hat G$, in such a way we obtain a functional, which is a generalization of $\Gamma$-function. It turns out that this functional, denoted as $\hat{\Gamma}_{\mu}$, which depends on some function $\mu$, such that $\mu(u)>0$ on $[0,2\pi]$, has a property which generalizes the famous property of the $\Gamma$-function: $\Gamma(x+1)=x\Gamma(x)$. Namely, we have:
\begin{eqnarray}
&&\int_0^{2\pi} g(v)\mu(v)\hat\Gamma_{\mu}(z+\delta_v)dv=\int _{0}^{2\pi}g(v)z(v)dv\hat\Gamma_{\mu}(z)+\nonumber\\
&&\frac{1}{t}\int_0^{2\pi}g''(v)\frac{\delta}{\delta z(v)}\hat\Gamma_{\mu}(z)dv, 
\end{eqnarray}
where $g(v)$ is any twice differentiable function on $[0,2\pi]$, such that $g(0)=g(2\pi)=0$ and 
$\delta_v=\delta(u-v)$ is a delta-function on the interval $[0,2\pi]$.

The outline of the paper is as follows. In Section 2 we remind the representation theory of $ax+b$ group and its relation to the theory of special functions. Section 3 contains basic facts about the classical Wiener measure. 
In Section 4 we define the path and loop groups  associated with $ax+b$ group and construct unitary 
representations of those. Section 5 is devoted to the relation of these representations to the 
loop generalization of the $\Gamma$-function. In Section 6, some of the many possible directions of further study are discussed.

\section{Unitary representations of the affine group of a line and $\Gamma$-function}
In this section we remind the reader of all the necessary facts about the unitary representations of $ax+b$ group, which is the affine group of the real line. In other words, each group element $g=g(a,b)$ is determined by the pair of real numbers $a,b$ such that $a>0$. The composition law is defined as follows: $g(a_1,b_1)g(a_2,b_2)=g(a_1a_2, a_1b_2+b_1)$. In the following we will call this group $G$.

The unitary representations of this group are constructed by means of method of induced representations. Consider the representation $r_{\lambda}$ of the subgroup $B$ of $G$ generated by $g(1,b)$-elements, such that $r_{\lambda}(g(1,b))=e^{\lambda b}$ (here $\lambda$ is complex). Then according to the construction of induced representations we have to consider the space of complex valued functions on $G$, i.e. $f(g(a,b))\equiv f(a,b)$ such that
\begin{eqnarray}
f(a,b+b_0)=e^{\lambda b_0} f(a,b).
\end{eqnarray}
Therefore $f(a,b)= e^{\lambda b}f(a,0)$ i.e. the function $f$ can be expressed in terms of the function $\phi(a)=f(a,0)$ on the subgroup $A$ generated by $g(a,0)$-elements. 
Then the operator of induced representation $R_{\lambda}$ acts on $f(g)$ as $R_{\lambda}(g_0)f(g)=f(gg_0)$, and the resulting formula for representation on the functions $\phi(a)$ is:
\begin{eqnarray}
R_{\lambda}(g_0)\phi(a)=e^{\lambda a b_0}\phi(a_0a),
\end{eqnarray}
or, in other words, since the space of functions $\phi$ is just a space of functions on the ray $0<x<\infty$, we have
$R_{\lambda}(g)\phi(x)=e^{\lambda bx}\phi(ax)$. One can consider the invariant measure on $A$, it is 
just $\frac{dx}{x}$. Then the following theorem is valid.\\

\noindent{\bf Theorem 2.1.} {\it i) Representation $R_{\lambda}$ of $G$ is unitary on $L^2(\mathbb{R}_{+},\frac{dx}{x})$ if $\lambda\in i\mathbb{R}$.\\
ii) Representations $R_{\lambda}$ and $R_{e^{\xi}\lambda}$ are equivalent if $\xi\in \mathbb{R}$. \\
iii) Consider a semigroup $G_+$(resp. $G_-$), consisting of such $g(a,b)$, that $b>0$ (resp. $b<0$). Then $R_{\lambda}$ with $\lambda\in \mathbb{C}$, such that $Re\lambda<0$ (resp. $Re\lambda>0$) is a representation for $G_+$(resp. $G_-$) on $L^2(\mathbb{R}_{+},\frac{dx}{x})$.}\\

As a consequence, one can see that there are only three classes of unequivalent unitary representations $R_{\lambda}$: $R_{\pm i}$ and $R_0$. The representations $R_{\pm i}$ are irreducible, while $R_0$ decomposes into the direct integral of 1-dimensional representations $T_{\rho}$, such that $T_{\rho}(g(a,b))=a^{i\rho}$. One can show that any other irreducible representation of $G$ falls into one of the classes $R_{\pm i}$,  $T_{\rho}$.

In order to generalize representations $R_{\lambda}$ to the loop case, it is useful to consider another form of the representations $R_{\lambda}$. Namely, instead of the space $L^2(\mathbb{R}_{+},\frac{dx}{x})$ one can 
consider $L^2(\mathbb{R},dt)$, by substitution $x=e^t$. Therefore, the formula for the representation $R_{\lambda}$ can be rewritten as follows:
\begin{eqnarray}\label{rl}
R_{\lambda}(g(e^{\alpha},b))f(t)=e^{i\lambda be^t}f(t+\alpha),
\end{eqnarray}
where $f \in L^2(\mathbb{R},dt)$ and we represented $a$ as $e^{\alpha}$. 

Now we discuss the relation of the representation $R_{\lambda}$ and $\Gamma$-function. It is well known, that a 
bilateral Laplace transform (or Fourier transform in the complex domain)
\begin{eqnarray}
\mathcal{L}f(p)=\frac{1}{\sqrt{2\pi}}
\int_{\mathbb{R}}e^{ipt}f(t)dt,
\end{eqnarray}
where $p$ is a complex number, has the inverse:
\begin{eqnarray}
\mathcal{L}^{-1}g(t)=\frac{1}{\sqrt{2\pi}}
\int_{\mathbb{R}+iT}e^{-ipt}g(p)dp,
\end{eqnarray}
where $T$ is a real number, so that the contour of integration is in the region of convergence of $g(p)$.
Therefore, if one can make sense of  $\mathcal{L}R_{\lambda}\mathcal{L}^{-1}$, it gives us a representation, equivalent 
to $R_{\lambda}$. Let $\mathcal{D}\subset L^2(\mathbb{R})$ be the space of $C^{\infty}$ functions with finite support. \\

\noindent{\bf Proposition 2.1.} \cite{vilenkin} {\it The bilateral Laplace transform of an element of $\mathcal{D}$ is an analytic function in the entire complex plane of exponential type, i.e. $|\mathcal{L}f(x+iy)|<Ce^{b|y|}$, where $C>0$, 
$b>0$. At the same time, the inverse Laplace transform of such analytic function belongs to $\mathcal{D}$. Moreover, we have the following property:
\begin{eqnarray}
\int_{\mathbb{R}}|\mathcal{L}f(x+iy)|^2dx<\infty,
\end{eqnarray}
where $f\in \mathcal{D}$.}\\

\noindent We notice that $\mathcal{D}$ is invariant under the action of the operators $R_{\lambda}(g)$. 
Let $\Gamma(z)$ denote the $\Gamma$-function: 
\begin{eqnarray}
\Gamma(z)=\int^{\infty}_0 e^{-x}x^{z-1}dx.
\end{eqnarray}
Then we have a Proposition.\\

\noindent {\bf Proposition 2.2.} \cite{vilenkin} {\it i) Consider the action of the $\mathcal{L}R_{\lambda}(g)\mathcal{L}^{-1}$ on $\mathcal{L}\mathcal{D}$, when $g=g(a,b)\in G_+$ and $\lambda<0$. Then 
\begin{eqnarray} 
\mathcal{L}R_{\lambda}(g)\mathcal{L}^{-1}f(t_1)=
\frac{1}{2\pi } \int_{\mathbb{R}+i0}\Gamma(it_1-it_2)a^{-it_1}\Big(\frac{-\lambda b}{a}\Big)^{it_2-it_1}f(t_2)dt_2,
\end{eqnarray}
where $f$ is an analytic function on $\mathbb{C}$ of the exponential type. \\
ii) One can analytically continue the expression above with respect to $-\lambda b$ to all complex plane except for 
negative real axis. }\\

\noindent We note, that we intentionally subdivided Proposition ii) into two statements. As we will see, in the loop case we will have i) only.

\section{Classical Wiener Measure: a reminder}
\noindent {\bf Notation.} In the following we will use several functional spaces, so let us fix notations. Let $\mu$ be 
the $\sigma$-additive measure on some space $X$. Then we denote as $L^2(X,d\mu;k)$ the Hilbert space 
of square-integrable complex-valued functions (here $k$ stands for $\mathbb{C}$ or $\mathbb{R}$). Since in the most of cases we will deal with $L^2(X,d\mu;\mathbb{C})$, we drop $\mathbb{C}$, i.e. $L^2(X,d\mu)\equiv L^2(X,d\mu;\mathbb{C})$. The space of real-valued continuous functions on the interval $[a,b]$ will be referred to as $C[a,b]$. Finally, the space of real-valued absolutely continuous functions on $[a,b]$, i.e. differentiable functions from $C[a,b]$, whose derivative belong to $L^2([a,b],dx;\mathbb{R})$, where $dx$ is the Lebesgue measure on $[a,b]$ will be denoted as $C'[a,b]$.\\

\noindent{\bf 3.1. Brief review of the abstract approach.} In this section we collect all the necessary facts about the Wiener measure. For a more detailed exposition of this subject one can consult \cite{quo}, \cite{hida}, \cite{gelfand}, \cite{lee}, \cite{frenkel}. 

Abstract Wiener measure is a Gaussian measure on a Banach space with certain properties. 
The construction is as follows. Start from Hilbert space $\mathcal{H}$ and consider a Gaussian measure associated with the unital operator and the norm $||\cdot||$ on the Hilbert space, such that heuristically measure can be represented as follows:
\begin{eqnarray} 
d\tilde{w}^t\sim e^{-\frac{||x||^2}{2t}}[dx].
\end{eqnarray}
However, this Gaussian measure is not $\sigma$-additive. In order to make it $\sigma$-additive, one  has to 
consider a weaker norm $|\cdot|$ with certain conditions on it with respect to measure $d\tilde{w}^t$. Then one can consider a completion $\mathcal{H}$ with respect to the norm $| \cdot |$. This will give a Banach space $\mathcal{B}$. Now the gaussian measure $d\tilde w^t$ can be extended to the Banach space, where it becomes 
$\sigma$-additive. 

The important example, which we will consider in the next subsection, is constructed as follows. Take the Hilbert space $C_0'[0,2\pi]$ of real valued absolutely continuous functions, such that $x(0)=0$ for any $x\in C_0'[0,2\pi]$ and 
$\int_0^{2\pi}(x'(u))^2du<\infty$. The inner product is given by: $\langle x_1,x_2\rangle=\int_0^{2\pi}x_1'(u)x_2'(u)du$. 
Then one can consider the weaker norm on this space: $|x|={\rm sup}_{u\in[0,2\pi]} x(u)$. The completion of $C_0'[0,2\pi]$ with respect to $|\cdot|$ is the Banach space $C_0[0,2\pi]$ of real valued continuous functions such that $x(0)=0$ for any $f\in C_0[0,2\pi]$. 

It appears that the norm $|\cdot|$ satisfies all necessary properties and therefore there exists a $\sigma$-additive Gaussian measure $dw^t$ on $C_0[0,2\pi]$. This measure is called classical Wiener measure. In the next subsection we will give its direct construction, using another approach.\\

\noindent {\bf 3.2. The construction of the classical Wiener measure and its basic properties.} Consider the space of continuous functions $C_0[0,2\pi]$. Let $C_{0,X}[0,2\pi]$ denote the closed subspace of  $C_0[0,2\pi]$, consisting 
of real valued continuous functions such that $x(2\pi)=X$ for some $X\in \mathbb{R}$. 

The following subsets of $C_0[0,2\pi]$ are called cylinder sets:
\begin{eqnarray} 
\{x\in C_0[0,2\pi]: x(\tau_1)\in A_1, \dots , x(\tau_n) \in A_n, 0<\tau_1<\dots\tau_n\le 2\pi\},
\end{eqnarray}
where $A_1,A_2,\dots A_n$ are Borel subsets of $\mathbb{R}$.\\

\noindent {\bf Definition 3.1.} {\it i)The Wiener measure with variance $t\ge 0$ is defined on the cylindrical sets of $C_0[0,2\pi]$ as follows:
\begin{eqnarray}
&&w^t( x(\tau_1)\in A_1, \dots , x(\tau_n) \in A_n)=\nonumber\\
&&\int_{A_1}\dots \int_{A_n} u_t(\Delta x_1, \Delta \tau_1)\dots 
u_t(\Delta x_n, \Delta \tau_n)dx_1\dots dx_n.
\end{eqnarray}
ii) The conditional Wiener measure of variance $t>0$ is defined on the cylinder sets of $C_{0,X}[0,2\pi]$ as follows:
\begin{eqnarray}
&&w^t_X( x(\tau_1)\in A_1, \dots , x(\tau_n) \in A_n)=\nonumber\\
&&\int_{A_1}\dots \int_{A_{n-1}} u_t(\Delta x_1, \Delta \tau_1)\dots 
u_t(\Delta x_n, \Delta \tau_n)dx_1\dots dx_{n-1},
\end{eqnarray}
where $u_t(x,s)=\frac{1}{2\pi\sqrt {ts}}e^{-\frac{x^2}{4\pi st}}$, $dx$ stands for Lebesgue measure on $\mathbb{R}$, 
$\Delta x_k=x_k-x_{k-1}$, $\Delta \tau_k=\tau_k-\tau_{k-1}$, $x_0=0$, $\tau_0=0$. In the case ii) $x_n=X$, $\tau_n=2\pi$.}\\

\noindent It appears that the resulting measures are $\sigma$-additive and moreover, we have the following theorem. \\

\noindent {\bf Theorem 3.1.} {\it i)Measure $w^t$ (resp. $w^t_X$) is $\sigma$-additive on the $\sigma$-field generated by the cylinder sets in $C_0[0,2\pi]$ (resp. $C_{0,X}[0,2\pi]$). The $\sigma$-field generated by the cylinder sets in $C_0[0,2\pi]$ (resp. $C_{0,X}([0,2\pi])$) is the Borel field of this Banach space.\\
ii) $w^t(C_0[0,2\pi])=1$, $w^t_X(C_{0,X}[0,2\pi])=u_t(X,2\pi)$.
}\\

One can show that the description of the Wiener measure that we gave in subsection 3.1. and the one introduced in this subsection agree (see Section I.5. of \cite{quo}). 
Now we want to relate $L^2$ spaces with respect to Wiener measure and conditional Wiener measure. In order to do that, notice the following property. If $f$ is an integrable function on $C_{0,X}([0,2\pi])$ for almost all $X$, then by the definition of the Wiener measure one has:
\begin{eqnarray}
\int_{C_{0}[0,2\pi]} f(x)dw^t(x)=\int_{\mathbb{R}}\Big(\int_{C_{0,X}[0,2\pi]}f(x)dw_X^t(x)\Big)dX.
\end{eqnarray}
Therefore, we have the proposition, which gives a decomposition of the $L^2$ space as a direct integral (for the 
definition of direct integral see e.g. \cite{vilenkin}, \cite{naimark}).\\

\noindent{\bf Proposition 3.1.}{\it The Hilbert space of square-integrable functions with respect to $dw^t$ decomposes as a direct integral:
\begin{eqnarray}
L^2(C_{0}[0,2\pi], dw^t)=\int_{\mathbb{R}}^{\oplus}L^2(C_{0,X}([0,2\pi]), dw^t_X)dX.
\end{eqnarray}
}
The following property, which will be crucial in this paper, describes the translation property of the Wiener measure and conditional Wiener measure.\\

\noindent {\bf Proposition 3.2.} {\it i) Let f be an integrable function on  $C_{0}([0,2\pi])$ and $y\in C'_0[0,2\pi]$, then 
\begin{eqnarray}
&&\int_{C_{0}[0,2\pi]} f(x)dw^t(x)=\nonumber\\
&&\int_{C_{0}[0,2\pi]} f(x+y)e^{-\frac{1}{t}\int^{2\pi}_0y'(u)dx(u)-\frac{1}{2t}\int^{2\pi}_0y'(u)y'(u)du}dw^t(x)
\end{eqnarray}
where $\int^{2\pi}_0y'(u)dx(u)$ is a Stiltjes integral.\\
ii) Let f be an integrable function in $C_{0,X}([0,2\pi])$ and $y$ as in i), such that $y(2\pi)=Y$ then,
\begin{eqnarray}
&&\int_{C_{0,X+Y}[0,2\pi]} f(x)dw_{X+Y}^t(x)=\nonumber\\
&&\int_{C_{0,X}[0,2\pi]} f(x+y)e^{-\frac{1}{t}\int^{2\pi}_0y'(u)dx(u)-\frac{1}{2t}\int^{2\pi}_0y'(u)y'(u)du}dw_{X}^t(x).
\end{eqnarray}
}
The space of Wiener measure translations, i.e. $C'_0[0,2\pi]$ is usually called Cameron-Martin space.\\

\noindent{\bf 3.3. Fourier-Wiener transform.} One can define a unitary operator on $L^2$ space for any abstract Wiener measure, which is similar to Fourier transform. We will need it only in the case of measure $dw^t_0$ on $C_{0,0}[2\pi]$.
The formula is as follows:
\begin{eqnarray}
Ff(y)=\int_{C_{0,0}[0,2\pi]}f(x+iy)dw^{2t}_0(x)
\end{eqnarray} 
Here $f\in L^2(C_{0,0}[0,2\pi], dw^t_0)$. The inverse transformation is given by the changing sign of $y$ in the formula 
above.\\

\noindent{\bf Theorem 3.2.} {\it The operator $F$ is unitary on $L^2(C_{0,0}[0,2\pi], dw^t_0)$.} 

\section{Unitary representations of the path and loop versions of $ax+b$ group}
\noindent{\bf 4.1. Path groups and loop groups associated to $G$.} In this subsection we will define all the path and loop groups, which we study in this article. 

At first we define the path group $PG$. Let us consider a set of elements of the form $g(e^{\alpha}, b)$, where $\alpha$ and $b$ are real valued absolutely continuous functions on $[0, 2\pi]$. We can define the bilinear operation as follows:
\begin{eqnarray}\label{mult}
g(e^{\alpha_1}, b_1)\cdot g(e^{\alpha_2}, b_2)=g(e^{\alpha_1+\alpha_2}, e^{\alpha_1}b_2+b_1),
\end{eqnarray}
where $e^{\alpha_1}b_2$ stands for pointwise multiplication of the absolutely continuous functions $e^{\alpha_1}$ and $b_2$. It is clear that the operation is well defined and satisfies the group laws. 
Let us denote the group which is a set of all such elements g($e^{\alpha}, b$) as $PG$. 

In the following we will sometimes need $\alpha$ to be in Cameron-Martin space for Wiener measure, so it is useful to introduce the based path group $P_0G\subset PG$, generated by $g(e^{\alpha},b)$, where $\alpha(0)=0$.

Similarly, one can define loop groups $\Omega G$ and $\Omega_0 G$. The loop group $\Omega G\subset PG$ is generated by $(e^{\alpha},b)$, where $\alpha, b$ are such that $\alpha(0)=\alpha(2\pi)$, $b(0)=b(2\pi)$, while 
the based loop group $\Omega_0 G\equiv P_0G\cap\Omega G$.

Finally, one can define the central extended versions of $\Omega G$ (resp. $\Omega_0 G$). Let us consider the elements $g(e^{\alpha}, b, s)$, where $s\in \mathbb{R}$ and $(e^{\alpha}, b)\in \Omega G$ (resp. $\Omega_0 G$). 
Then the multiplication law can be modified in the following way:
\begin{eqnarray}
&&g(e^{\alpha_1}, b_1, s)\cdot g(e^{\alpha_2}, b_2, t)=\nonumber\\
&&g(e^{\alpha_1+\alpha_2}, e^{\alpha_1}b_2+b_1, t+s+
k\int_0^{2\pi}\alpha_1(u)\alpha'_2(u)\ud u),
\end{eqnarray}
where $k\in \mathbb{R}$. We denote the corresponding group $\hat{G}$ (resp. $\hat{G}_{0}$) and $k$ is called central charge.\\

\noindent{\bf 4.2. Unitary representations.} 
In this subsection we will describe the unitary representations 
of both $PG$ and $\Omega G$, as well as their subgroups $P_0G$ and $\Omega_0 G$, which are the appropriate generalizations of the unitary representations of $G$ which we considered in Section 2. 

We will start from $P_0G$ and use the same approach as before, i.e. we will use the method of induced representations. Let us take one-dimensional representations of the $B$-subgroup (i.e. subgroup of elements of the 
form $g(1,b)$) of the form $\tilde{r}_{\lambda}(g(1,b))=e^{\int_0^{2\pi}\lambda(u)b(u)\ud u}$, where $\lambda\in L^2[0,2\pi]$. Following the method of induced representations, like we did in the case of group $G$, we arrive to the following formula for the representations on 
the space of functionals on the $A$-subgroup:
\begin{eqnarray}
\tilde{R}_{\lambda}(g(e^{{\alpha}_0},{b_0}))f(\alpha)=e^{\int_0^{2\pi}\lambda(u)b_0(u)e^{{\alpha}(u)}\ud u}f(\alpha+\alpha_0).
\end{eqnarray}
However, in order to make these representations unitary one needs to define the proper inner product on the space 
of functionals. In order to do that, one has to consider Wiener measure. It is not invariant under translations, so one should improve the formula for the representations in order to make them unitary.  
Let us consider the Hilbert space $H^0_p=L^2(C_0[0,2\pi], dw^t;\mathbb{C})$. The following statement is true.\\

\noindent {\bf Theorem 4.1.} {\it The following action of $P_0G$ on the space of functionals of continuous functions  
\begin{eqnarray}\label{act}
&&\rho_{\lambda}(g(e^{\alpha},b))(f)(x)=\nonumber\\
&&e^{-\frac{1}{4t}\int_0^{2\pi}\alpha'(u)\alpha'(u)\ud u-\frac{1}{2t} \int_0^{2\pi}\alpha'(u)\ud x(u)}
e^{\int_0^{2\pi} \lambda(u)b(u)e^{x(u)}\ud u}f(x+\alpha)
\end{eqnarray}
defines the unitary representation of $P_0G$ on $H^0_p$ iff $i\lambda(u)\in L^2([0,2\pi];\mathbb{R})$.}\\

\noindent {\bf Proof.} First we prove that this action is actually an action of a group, i.e. 
$\rho_{\lambda}(g(e^{\alpha_1},b_1))\rho_{\lambda}(g(e^{\alpha_2},b_2))=\rho_{\lambda}(g(e^{\alpha_1+\alpha_2}, e^{\alpha_1}b_2+b_1))$. We 
check it, writing the explicit expressions:
\begin{eqnarray}
&&\rho_{\lambda}(g(e^{\alpha_1},b_1))\rho_{\lambda}(g(e^{\alpha_2},b_2))f(x)=\nonumber\\
&&\rho_{\lambda}(g(e^{\alpha_1},b_1))(e^{-\frac{1}{4t}\int_0^{2\pi}\alpha_2'(u)\alpha_2'(u)\ud u-\frac{1}{2t} \int_0^{2\pi}\alpha_2'(u)\ud x(u)}\nonumber\\
&&e^{\int_0^{2\pi} \lambda(u)b_2(u)e^{x(u)}\ud u}f(x+\alpha_2))=\nonumber\\
&&e^{-\frac{1}{4t}\int_0^{2\pi}\alpha_2'(u)\alpha_2'(u)\ud u-\frac{1}{4t}\int_0^{2\pi}\alpha_1'(u)\alpha_1'(u)\ud u-\frac{1}{2t}\int_0^{2\pi}\alpha_1'(u)\alpha_1'(u)\ud u}\nonumber\\
&&e^{-\frac{1}{2t} \int_0^{2\pi}(\alpha_1'(u)+\alpha_2'(u))\ud x(u)}e^{\int_0^{2\pi} \lambda(u)(b_2(u)e^{\alpha_1(u)}+b_1(u))e^{x(u)}\ud u}f(x+\alpha_1+\alpha_2))=\nonumber\\
&&\rho_{\lambda}(e^{\alpha_1+\alpha_2}, e^{\alpha_1}b_2+b_1)f(x).
\end{eqnarray}
Now we prove that the action $\rho_{\lambda}$  defines a unitary representation for $P_0G$. In order to do that we just need to use the definition of unitarity and the translation invariance property of the Wiener measure:
\begin{eqnarray}
&& (\rho_{\lambda}(g(e^{\alpha},b))(f),\rho_{\lambda}(g(e^{\alpha},b)) (h))_{H_p}=\nonumber\\
&&\int_{C_0[0,2\pi]} \overline{\rho_{\lambda}(g(e^{\alpha},b))(f)(x)}\rho_{\lambda}(g(e^{\alpha},b)) (h)(x)dw^t(x)=\nonumber\\
&&\int_{C_0[0,2\pi]}  e^{-\frac{1}{2t}\int_0^{2\pi}\alpha'(u)\alpha'(u)\ud u-\frac{1}{t}\int_0^{2\pi}\alpha'(u)\ud x(u)}\overline{f(x+\alpha)}g(x+\alpha)dw^t(x)=\nonumber\\
&&\int \overline{f(x)}g(x)dw(x)=(f,g)_{H_p}.
\end{eqnarray}
Thus the Theorem is proven.\hfill$\blacksquare$\\

\noindent In the following we will denote representations $\rho_{\lambda}(g(e^{\alpha},b))\equiv \rho_{\lambda}(e^{\alpha},b)$. 
 
 Using the same steps as in the proof of the Theorem above one can show that analogous fact for the Hilbert space $H^{0,X}_p=L^2(C_{0,X}[0,2\pi], dw_X^t;\mathbb{C})$ and group $\Omega_0 G$ is true.\\

\noindent {\bf Proposition 4.1.} {\it i) Formula \rf{act} defines a unitary representation $\rho^X_{\lambda}$ of $\Omega_0 G$ on $H^{0,X}_p$.\\
ii) Representation $\rho^{X_1}_{\lambda}$ on $H^{0,X_1}_p$ is equivalent to the representation of $\rho^{X_2}_{e^{\eta}\lambda}$ on $H^{0,X_2}_p$, where $\eta$ is absolutely continuous, $\eta(0)=0$, $\eta(2\pi)=X_1-X_2$.}\\

\noindent {\bf Proof.} The first part follows from the proof of the Theorem above, while the second part is the direct consequence of the translation formula.\hfill$\blacksquare$\\

\noindent Therefore, all representations $\rho^X_{\lambda}$ are equivalent to the representations $\rho^0_{e^{\xi}\lambda}$ on 
$H^0_l\equiv L^2(C_{0,0}[0,2\pi], dw^t;\mathbb{C})$, where $\xi$ is absolutely continuous, $\xi(0)=0, \xi(2\pi)=X$. 
Moreover, one can show that this equivalence does not depend on the choice of such $\xi$. 

However, we note here, that since $\Omega_0 G\subset P_0G$, the Hilbert space $H^{0}_p$ is a representation space for $\Omega_0 G$ too. Therefore we have a Proposition.\\

\noindent {\bf Proposition 4.2.} {\it The representation $\rho_{\lambda}$ of $\Omega_0G$ on $H^0_p$ is the direct integral of the representations of $\Omega_0G$: 
\begin{eqnarray}
\rho_{\lambda}=\int^{\oplus}_{\mathbb{R}} \rho^{X}_{\lambda}dX, 
\end{eqnarray}
where $\rho^{X}_{\lambda}$ is equivalent to $\rho^0_{\lambda_X}$ such that  
$\lambda_X=e^{\xi}\lambda$, $\xi\in C'[0,2\pi]$ and $\xi(0)=0$, $\xi(2\pi)=X$.\\
}

\noindent {\bf Proof.} 
The proof follows from the decomposition of  $H^0_p$:
\begin{eqnarray}
H^0_p=\int^{\oplus}_{\mathbb{R}} H^{0,X}_pdX
\end{eqnarray}
and Proposition above. \hfill $\blacksquare$\\

\noindent In order to define the representation of $PG$ and $\Omega G$ one needs to extend the Hilbert space, i.e. one has to consider the Hilbert spaces $H_p=L^2(C_{0}[0,2\pi]\oplus \mathbb{R}, dw^t\times dx_0;\mathbb{C})$ and 
$H_l=L^2(C_{0,0}[0,2\pi]\oplus \mathbb{R}, dw^t\times dx_0;\mathbb{C})$, where $dx_0$ is the Lebesgue measure on $\mathbb{R}$.
We also note that every element $(e^{\alpha}, b)$ of $PG$ (resp. $\Omega G$) can be represented as 
$(e^{\alpha}, b)=(e^{\alpha_0}e^{\t \alpha}, b)$, where $\alpha_0=\alpha(0)$ and $\t \alpha(0)=0$.  
Then the following theorem is true.\\

\noindent {\bf Theorem 4.2.} {\it i) The following action of the group element $g(e^{\alpha}, b)$ on the space of 
functionals of $x(u), x_0$:
\begin{eqnarray}\label{act2}
&&\t{\rho^p}_{\lambda}(e^{\alpha},b)(f)(x, x_0)=\\
&&e^{-\frac{1}{4t}\int_0^{2\pi}\alpha'(u)\alpha'(u)\ud u-\frac{1}{2t} \int_0^{2\pi}\alpha'(u)\ud x(u)}
e^{\int_0^{2\pi} \lambda(u)b(u)e^{x(u)+x_0}\ud u}f(x+\t\alpha, x_0+\alpha_0)\nonumber
\end{eqnarray}
defines a unitary representation of $PG$  on the Hilbert space $H_p$ iff $\lambda\in iL^2([0,2\pi]$. 
ii) The formula \rf{act2} defines the unitary representation 
$\t \rho^l_{\lambda}$ of $\Omega G$ on the Hilbert space $H_l$ iff $\lambda \in iL^2([0,2\pi];\mathbb{R})$.}\\

\noindent Note that representations $\rho^p_{\lambda}$,  $\rho^l_{\lambda}$ are equivalent to 
$\rho^p_{e^t\lambda}$,  $\rho^l_{e^t\lambda}$ correspondingly, for any $t\in \mathbb{R}$. 
In particular, in the case when $\lambda=const$ 
there exist only three inequivalent classes of unitary representations of $PG$ ($\Omega G$), 
corresponding to  $\t \rho_0$, $\t\rho_i$, $\t\rho_{-i}$ respectively, like it was in the finite dimensional case. \\

\noindent{\bf 4.3. Central extension and Lie algebra generators.} In this subsection we will construct the unitary representations of groups 
$\hat{G}$ on $H_l$. Let $(e^{\alpha},b, s)$ be an element of $\hat G$. Let us consider the following action of this element on the function $H_l$:
\begin{eqnarray}\label{act3}
&&\rho^l_{\lambda,k}((e^{\alpha},b, s))(f)(x,x_0)=
e^{is}e^{-\frac{1}{4t}\int_0^{2\pi}\alpha'(u)\alpha'(u)\ud u-\frac{1}{2t} \int_0^{2\pi}\alpha'(u)\ud x(u)}\nonumber\\
&&e^{ik\int_0^{2\pi}\alpha(u)dx(u)}e^{\int_0^{2\pi} \lambda(u)b(u)e^{x(u)+x_0}\ud u}f(x+\t \alpha,x_0+\alpha_0).
\end{eqnarray}
{\bf Theorem 4.3.} {\it Formula \rf{act3} defines a unitary representation of $\hat G$ for the central charge $k\in \mathbb{R}$, and $\lambda\in i L^2([0,2\pi];\mathbb{R})$. Representation ${\rho^l}_{e^{\xi}\lambda,k}$ is equivalent to ${\rho^l}_{\lambda,k}$ if $\xi$ is any absolutely continuous function on $[0,2\pi]$, such that $\xi(0)=\xi(2\pi).$}\\ 

Ignoring $x_0$ dependence in \rf{act3}, one can define the unitary representation of $\hat{G}_0$ on $H^0_l$, for 
which we will keep the same notation $\rho^l_{\lambda,k}$.  

\noindent We also mention the following important statement.\\

\noindent {\bf Proposition 4.3.} {\it Operators $\rho^l_{\lambda,k}(e^{\alpha_1}, 0, s_1)$ and $\rho^l_{\lambda,-k}(e^{\alpha_2}, 0, s_2)$  corresponding to the action of the $A$ subgroup of $\hat{G}$ commute.}\\

\noindent The proof follows from the formula \rf{act3} by direct calculation. \\ 
One can write down the expressions for the Lie algebra generators. As we can see, the generator which is responsible for the action of the "$B$-subgroup"  has the form:
\begin{eqnarray}
T_bf(x,x_0)=\frac{d}{d\epsilon}_{\vert_{\epsilon=0}}\rho^l_{\lambda,k}(1,\epsilon b,0)=\int_0^{2\pi}\lambda(u)b(u)e^{x(u)+x_0}du.
\end{eqnarray}
It is well-defined for all the functionals $f(x,x_0)\in H_l$ such that $e^{x_0}f(x,x_0)\in H_l$. However, the $a$-subgroup generators are more peculiar:
\begin{eqnarray}
&&D_{k,\alpha}f(x)=\frac{d}{ds}_{\vert_{\epsilon=0}}\rho^l_{\lambda,k}(e^{\epsilon\alpha},0,0)=\\
&&(ik\int_0^{2\pi}\alpha(u) dx(u)+\frac{1}{2t}
\int_0^{2\pi}\alpha'(u)dx(u))f(x)+\frac{d}{d\epsilon}_{|_{\epsilon=0}}f(x+\epsilon\t \alpha, x_0+\epsilon \alpha_0).\nonumber
\end{eqnarray}
Therefore the action of the Lie algebra element $T_{\alpha}$ is defined only on the set of weakly differentiable functionals, which form a dense subset in $H_l$. In the notations of variational calculus, one can 
introduce the operators:
\begin{eqnarray}\label{var}
&&D_{\alpha,k}(\alpha)=\alpha_0\frac{\p}{\p x_0}+\int_0^{2\pi}\t\alpha(u)\frac{\delta}{\delta x(u)}du+
\frac{1}{2t}\int_0^{2\pi}\t\alpha'(u)x'(u)du+\nonumber\\
&& ik\int_0^{2\pi} du\t\alpha(u)x'(u)du,  \nonumber\\
&&T_b=\int_0^{2\pi}\lambda(u)b(u)e^{x_0+x(u)}du.
\end{eqnarray}
 One can see that \rf{var} satisfies the needed relation:
\begin{eqnarray}
[D_{\alpha_1,k}, D_{\alpha_2,k}]=-2ik \int_0^{2\pi}\alpha_1'(u)\alpha_2(u)du, 
\quad [D_{\alpha,k}, T_b]=T_{\alpha b}.
\end{eqnarray}
Finally, we notice that 
\begin{eqnarray}
[D_{k,\alpha_1}, D_{-k,\alpha_2}]=0.
\end{eqnarray}
This is just a paraphrasing of the Proposition above on the Lie algebra level.\\
  
Finally, let us consider the following two semigroups of $\hat{G}_+,\hat{G}_-\in\hat G$, such that they 
consist of group elements $g(e^\alpha, b)$, where $b(u)>0$ or $b(u)<0$ for all $u$ correspondingly. We mention the following important statement, which we will use in the next section.\\
 
\noindent{\bf Proposition 4.4.} {\it Let $Im\lambda, Re\lambda\in L^2([0,2\pi];\mathbb{R})$. Then \rf{act3} defines a representation of 
the semigroup $\hat{G}_+$ (resp. $\hat{G}_-$) if $Re \lambda(u)<0$ (resp.  $\lambda>0$) on $[0,2\pi]$.}\\
  
\noindent In the following we will denote the representations of $\hat{G}_+,\hat{G}_-$ the same way, namely 
$\rho^l_{\lambda,k}$.\\ 

\noindent{\bf 4.4. Irreducibility and classification of unitary representations.} In Section 2 we learned that the 
irreducible unitary representations of $G$ turn out to be equivalent to either $R_{\pm i}$ or $T_{\rho}$. In this subsection we will show that representations $\rho^l_{\lambda,k}$ such that $\lambda(u)\equiv \lambda$, where $\lambda\in \mathbb{R}\backslash 0$, are irreducible. Moreover, we will talk about the classification of all irrieducible unitary representations of $\hat{G}$. \\

\noindent{\bf Theorem 4.4.} {\it Representations $\rho^l_{\lambda,k}$, such that  $\lambda(u)\equiv\lambda\in \mathbb{R}\backslash 0$ of $\hat{G}$ ($\hat{G}_0$) on $H_l$ ($H^0_l$), are irreducible.}\\

\noindent{\bf Proof.} We will prove this fact for $\hat{G}_0$, because the proof for $\hat{G}$ requires just a minor modification. 

If $\rho^l_{\lambda,k}$ is not irreducible then there exists such bounded operator $P$, which commutes with the action of $\hat{G}_0$. Let us show that such the operator $P$ should be proportional to identity operator. 
Due to the continuity of $P$ we have the relation $P\mathcal{R}_b=\mathcal{R}_bP$, where $\mathcal{R}_b$ is a bounded operator on $H^0_l$ corresponding to the multiplication of the functional $e^{i\int_0^{2\pi}b(u)e^{x(u)}du}\in H^0_l$ for any $b(u)\in L^2([0,2\pi];\mathbb{R})$. 
Let us consider the expression:
\begin{eqnarray}\label{seq}
f_N(x)=\int_{C_{0,0}[0,2\pi]}\sum^N_{n=1}c_n e^{\int_0^{2\pi}\xi_n(u)y(u)du+i\int_0^{2\pi}e^{x(u)}\xi_n(u)du} dw^0_{2t}(y),
\end{eqnarray}
where $\xi_n$ are any continuous functions on $[0,2\pi]$ and $c_n$ are some constants. 
The function $f_N(x)$ is bounded and commutes with $P$. One can see that the bounded functional $g_N(s)$ such that 
$f_N(x)=g_N(e^x)$ is actually a Fourier-Wiener transform of the function $h_N(y)=\sum^N_{n=1}c_n \exp(\int_0^{2\pi}\xi_n(u)y_n(u)du)$. Now consider any functional $g(x)\in H_l^0$, which is bounded, then 
$g(e^x)=f(x)$ is also bounded. Therefore we see that one can construct any functional of the form  $g(e^x)$  on $C_{0,0}[0,2\pi]$ as a limit of certain functionals $f_N$ like in \rf{seq}, because the space spanned by exponential functions $e^{\int_0^{2\pi}\xi_n(u)y(u)du}$, where $\xi_n\in L^2([0,2\pi],\mathbb{R})$, is a dense set in $H_l^0$.  
Since $P$ is continuous, this means that it commutes with any bounded functional of the form $g(e^x)$. At the same time, these functionals are dense in the set of all bounded functionals $f(x)$. In order to see that, take a functional 
$s_n=e^{\frac{i}{n}\int_0^{2\pi}\eta(u)(e^{nx(u)}-1)du}$. Then $\lim_{n\to \infty}s_n=e^{i\int_0^{2\pi}\eta(u)x(u)du}$. Since the space spanned by exponential functionals of the form $e^{i\int_0^{2\pi}\eta(u)x(u)du}$, where $\eta\in L^2([0,2\pi],\mathbb{R})$ is dense in $H_l^0$, then bounded functionals of the form  $g(e^x)$ are dense in the space of all bounded functionals of $x$. Hence, $P$ has to commute with all bounded functionals and therefore, the action of $P$ on $H_l^0$ is equivalent to the multiplication on a function $a(x)$, which has to be bounded (since $P$ is bounded) 
a.e. with respect to $dw^0_t$. Since $P$ has to commute with the action of the $A$-subgroup we find that $a(x)$ 
can only be a constant. In such a way $\rho^l_{\lambda,k}$ is irreducible.
\hfill $\blacksquare$\\

One can prove that if $i\lambda$ is either a positive or negative function from $L^2([0,2\pi],\mathbb{R})$, 
the representation $\rho^l_{\lambda,k}$ is also irreducible. 
However, if $\lambda(u)=0$  $\rho^l_{\lambda,k}$ reduces to the representation of the $A$-subgroup, i.e. the representations of the loop Heisenberg group. It is not irreducible since $[\rho^l_{\lambda,k}(e^{\alpha_1}, 0, s), \rho^l_{\lambda,-k}(e^{\alpha_2}, 0, t)]=0$. To author's knowledge the classification of the unitary representations of the loop Heisenberg group is not yet known (see e.g. \cite{ccr} for review of the subject). 
The same argument as in the Theorem 4.4, shows that if $\lambda(u)=0$ for $u\in [a,b]$ the representation $\rho^l_{\lambda,k}$ is not irreducible. Therefore, we make the following conjecture about the classification of the irreducible representations 
of $\hat G$. 
Hence the irreducible unitary representations of $\hat G$ which we considered, are either equivalent to the  representations of the $A$-subgroup, or equivalent to the representations $\rho^l_{\lambda,k}$, where $i\lambda\in L^2([0,2\pi], \mathbb{R})$ is strictly positive or negative function on $[0,2\pi]$. Moreover, we know that  
$\rho^l_{\lambda,k}$ and $\rho^l_{e^{\xi}\lambda,k}$ are equivalent to each other if $\xi\in C'_{0,0}[0,2\pi]$. 
Therefore, an interesting problem to study is the classification of all finite-dimensional representations of $\hat{G}$. A reasonable conjecture (in the analogy with the finite-dimensional case) would be that the three discussed classes of representations, namely $\rho^l_{\lambda,k}$ for positive and negative $\lambda$ and irreducible unitary representations of loop Heisenberg group exhaust all irreducible unitary representations of $\hat{G}$. \\
 
\noindent {\bf Remark.} We hope (see the last section for more details) that the representations $\rho^l_{\lambda,k}$ of $\hat G$ for constant $\lambda$ are closed under suitable ''fusion'' tensor product in the analogy with representations $R_{\lambda}$ of $G$.

\section{(Loop) $\Gamma$-function and the action of the affine loop group}

\noindent{\bf 5.1. Fourier transform for the classical Wiener measure.} In this section we consider the generalizations of the formula relating the action of the group $G$ and the $\Gamma$-function (see section 1). In particular, we will introduce a new object, which we will refer to as $loop$ $\Gamma$-$function$. In order to do that one needs to construct the generalization of the Fourier/Laplace transform in the case of Wiener measure. We already have seen a  unitary transformation on the $L^2$ space for the abstract Wiener measure, called Fourier-Wiener transform, 
but we will choose another transformation here.

Let us consider the following transformation on the Hilbert space $H^0_l$:
\begin{eqnarray}
\mathcal{F}f(p)=\int_{C_{0,0}[0,2\pi]} e^{i\int_0^{2\pi}p(u)x(u)du} f(x)dw_0^t(x),
\end{eqnarray}
where $p(u)\in C[0,2\pi]$ and $p(0)=p(2\pi)=0$. Unlike the usual Fourier transform for Lebesgue measure on a real 
line, transformation $\mathcal{F}$ is not a unitary operator.\\

\noindent{\bf Proposition 5.1.} {\it The operator  $\mathcal{F}$ is a compact normal operator on  
$H^0_l$ with no zero eigenvalues.}\\

\noindent{\bf Proof.} The general condition \cite{conway} for the general integral operator\\  
$K: L^2(X,d\mu(x))\to L^2(Y,d\nu(y))$, such that 
\begin{eqnarray}
Kf(y)=\int_X K(y,x)d\mu(x)
\end{eqnarray}
to be compact is that $\int \int  |K(x,y)|^2d\mu(x)d\nu(y)<\infty$. In order to prove that it is normal, i.e. $\mathcal{F}\mathcal{F}^*=\mathcal{F}^*\mathcal{F}$ one just needs to write explicitly the resulting expressions and then use 
the Fubini theorem. To show that the operator $\mathcal{F}$ has no nonzero eigenvalues, one needs to use the fact that the exponentials of the form $e^{\int_0^{2\pi}{\alpha'(u)dx(u)}}$, where $\alpha(u)$ is absolutely continuous, form a dense subset in $H^0_l$ \cite{hida}.   
\hfill $\blacksquare$.\\

\noindent Taking into account that $\mathcal{F}$ is a normal operator and using the polar decomposition theorem, one can decompose it 
as $\mathcal{F}=U_{\mathcal{F}}K$, where $U_{\mathcal{F}}$ is a unitary operator on $H^0_l$ and $K$ is compact self-adjoint operator, such that $K=\sqrt{\mathcal{F}\mathcal{F}^*}$. 

It is obvious that the operators $\mathcal{F}$ and $U_{\mathcal F}$ can be continued to the space $H_l$. 
At the same time, we define another unitary operator $F$, which is a standard Fourier transform with respect to measure $dx_0$. One can show that $\mathcal{F}$ and $F$ commute. \\

\noindent{\bf 5.2. Loop $\Gamma$-function}.
Let us consider the  following expression:
\begin{eqnarray}
\rho^{l,\mathcal{F}}_{\lambda,k}(e^{\alpha},b,s)=U_{\mathcal{F}}\rho^l_{\lambda,k}(e^{\alpha},b,s)U^*_{\mathcal{F}}.
\end{eqnarray}
Since $U_{\mathcal{F}}$ is unitary, $\rho^{l,\mathcal{F}}_{\lambda,k}$ defines an equivalent representation 
of $\hat G$ on $H_l$. Because of the results of subsection 5.1, on the image  of $K$, one can rewrite it as follows:
\begin{eqnarray}
\rho^{l,\mathcal{F}}_{\lambda,k}(e^{\alpha},b,s)=K^{-1}\mathcal{F}\rho^l_{\lambda,k}(e^{\alpha},b,s)\mathcal{F}^*K^{-1}.
\end{eqnarray}
Here $K$ is a fixed self-adjoint operator, so we are interested in the object $\mathcal{F}\rho^l_{\lambda,k}(e^{\alpha},b)\mathcal{F}^*$. We consider the case when $(e^{\alpha},b)$ in $\hat{G}_{+}$ and $Im\lambda=0$, $Re \lambda(u)>0$. Let us write it down explicitly:
\begin{eqnarray}
&&\mathcal{F}\rho^l_{\lambda,k}(e^{\alpha},b,s)\mathcal{F}^*f(x,x_0)=\nonumber\\
&&e^{is}e^{-\frac{1}{4t}\int_0^{2\pi}\alpha'(u)\alpha'(u)du}\int_{C_{0,0}[0,2\pi]} e^{-i\int_0^{2\pi}(p(u)x(u))du}
e^{-\frac{1}{2t} \int_0^{2\pi}\alpha'(u)\ud p(u)}\nonumber\\
&&e^{ik\int_0^{2\pi}\alpha(u)dp(u)}e^{\int_0^{2\pi} \lambda(u)b(u)e^{p(u)+x_0}\ud u}\nonumber\\
&&\int_{C_{0,0}[0,2\pi]} e^{i\int_0^{2\pi}(p(u)+\t\alpha(u)y(u)du} f(y,x_0+\alpha_0)dw_0^t(y)dw_0^t(p).
\end{eqnarray}
Using the Fubini theorem, one can rewrite it as follows:
\begin{eqnarray}
&&\mathcal{F}\rho^l_{\lambda,k}(e^{\alpha},b,s)\mathcal{F}^*f(x,x_0)=\nonumber\\
&&\int_{C_{0,0}[0,2\pi]} \mathbb{K}^{\lambda,k}_{s,\alpha,b}(x-y, x_0) e^{i\int_0^{2\pi}\t\alpha(u)y(u)du}f(y, x+\alpha_0)dw_0^t(y),
\end{eqnarray}
where 
\begin{eqnarray}
&&\mathbb{K}^{\lambda,k}_{s,\alpha,b}(x-y,x_0)=e^{is}e^{-\frac{1}{4t}\int_0^{2\pi}\alpha'(u)\alpha'(u)}\\
&&
\int_{C_{0,0}[0,2\pi]} e^{i\int_0^{2\pi}p(u)(x(u)-y(u))du}e^{\int_0^{2\pi} \lambda(u)b(u)e^{p(u)+x_0}\ud u}\nonumber\\
&&e^{ik\int_0^{2\pi}\alpha(u)dp(u)}e^{-\frac{1}{2t} \int_0^{2\pi}\alpha'(u)\ud p(u)}
dw_0^t(p).\nonumber
\end{eqnarray}
In the case when $\alpha$ is twice differentiable, one can see that the object $\mathbb{K}^{\lambda,k}_{s,\alpha,b}(z,x_0)$ up to factors independent of $x,y$ is a particular case of the following functional:
\begin{eqnarray}
\hat\Gamma_{\mu}(z)=\int_{C_{0,0}[0,2\pi]} e^{\int_0^{2\pi}p(u)z(u)du}e^{-\int_0^{2\pi} \mu(u)e^{p(u)}\ud u}dw_0^t(p),
\end{eqnarray}
where $Re z, Im z \in L^2([0,2\pi];\mathbb{R})$, $\mu\in L^2([0,2\pi];\mathbb{R})$. We will call $\hat\Gamma_{\mu}(z)$ a {\it loop} {\it Gamma} $function$ or simply $\hat \Gamma$-$functional$. It has the following properties.\\

\noindent {\bf Theorem 5.1.} \\
\noindent{\it i) $\hat\Gamma_{\mu}(z)$ is well defined for any $Re z, Im z \in L^2([0,2\pi];\mathbb{R})$, $\mu\in L^2([0,2\pi];\mathbb{R})$ and $\mu(u)\ge 0$ on $[0,2\pi]$. \\
ii) The following relation is valid: 
\begin{eqnarray}\label{gamma}
&&\int_0^{2\pi} g(v)\mu(v)\hat\Gamma_{\mu}(z+\delta_v)dv=\int _{0}^{2\pi}g(v)z(v)dv\hat\Gamma_{\mu}(z)+\nonumber\\
&&\frac{1}{t}\int_0^{2\pi}g''(v)\frac{\delta}{\delta z(v)}\hat\Gamma_{\mu}(z)dv, 
\end{eqnarray}
where $g(v)$ is any twice differentiable function on $[0,2\pi]$, such that $g(0)=g(2\pi)=0$, 
$\delta_v=\delta(u-v)$ is a delta-function on the interval $[0,2\pi]$ and \\
$
\int_0^{2\pi} \xi(v)\frac{\delta}{\delta z(v)}\hat\Gamma_{\mu}(z)=\frac{d}{d\epsilon}_{|\epsilon=0}\hat\Gamma_{\mu}(z+\epsilon\xi)$ for any continuous function $\xi$. }\\

\noindent{\bf Proof.} An important step in the proof is the 
consideration of the infinitesimal form of the translation invariance. Let $f\in H^0_l$ such that it is weakly differentiable. Then we have the following property:
\begin{eqnarray}
&&\int_{C_{0,0}[0,2\pi]} e^{-\frac{1}{2t}\int_0^{2\pi}\epsilon^2g'(u)g'(u)\ud u+\frac{1}{t} \int_0^{2\pi}\epsilon g''(u)x(u)\ud u}f(x+\epsilon g )dw_0^t(x)=\nonumber\\
&&\int f(x)  dw^t(x),
\end{eqnarray}
where $g\in C_{0,0}^2[0,2\pi]$.
This is just the translation property (the translation is with respect to function $\epsilon g$), where $\epsilon$ is some real parameter.  Then, if we differentiate with respect to $\epsilon$ at zero, we obtain the following formula:
\begin{eqnarray}\label{rel}
&&\int_{C_{0,0}[0,2\pi]}\Big(\frac{1}{t} \int_0^{2\pi}g''(u)x(u)\ud u\Big)f(x)dw_0^t(x)=\nonumber\\
&&-\int_{C_{0,0}[0,2\pi]} \frac{d}{d\epsilon}_{|_{\epsilon=0}}f(x+\epsilon g)dw_0^t(x).
\end{eqnarray}
Let us apply this property to the integrand of the $\hat \Gamma$-functional, i.e. let 
$f(x)=F(x,z)\equiv e^{\int_0^{2\pi}x(u)z(u)du}e^{-\int_0^{2\pi} \mu(u)e^{x(u)}\ud u}$. Then we have:
\begin{eqnarray}\label{rel1}
&&\int_{C_{0,0}[0,2\pi]}\Big(\frac{1}{t} \int_0^{2\pi}g''(u)x(u)\ud u\Big)F(x,z)dw_0^t(x)=\nonumber\\
&&\frac{1}{t}\int_0^{2\pi}g''(v)\frac{\delta}{\delta z(v)}\hat\Gamma_{\mu}(z),
\end{eqnarray}
and
\begin{eqnarray}\label{rel2}
&&-\int_{C_{0,0}[0,2\pi]} \frac{d}{d\epsilon}_{|_{\epsilon=0}}F(x+\epsilon g,z)dw_0^t(x)=\nonumber\\
&&-\int_{C_{0,0}[0,2\pi]}\Big(\int^{2\pi}_0 g(v)z(v)dv-\int^{2\pi}_0\mu(v)g(v)e^{x(v)}dv\Big)F(x,z)dw_0^t(x)=\nonumber\\
&&-\int^{2\pi}_0 g(v)z(v)dv\hat\Gamma_{\mu}(z)+\int^{2\pi}_0g(v)\mu(v)\hat{\Gamma}(z+\delta_v)dv.
\end{eqnarray}
Therefore, combining \rf{rel}, \rf{rel1} and \rf{rel2} we obtain \rf{gamma}. 
\hfill $\blacksquare$\\

The property ii) from the theorem above is a natural generalization of the property of the ordinary $\Gamma$-function: 
$\Gamma(z+1)=z\Gamma(z)$. We notice, however, that there is also an extra term, depending on the $t$-parameter, related to the choice of the measure. So, the proper analogue of the functional $\hat{\Gamma}$ would be a "regularized" $\Gamma$-function:
\begin{eqnarray}
\Gamma_{\mu,t}(z)=\int_{\mathbb{R}} e^{-\mu e^x} e^{zx}e^{-\frac{x^2}{2t}}dx,
\end{eqnarray}
where $\mu, t\in{R_{+}}$. It satisfies the equation
\begin{eqnarray}
\mu\Gamma_{\mu,t}(z+1)=z\Gamma_{\mu,t}(z)-\frac{1}{t}\frac{d\Gamma(z)}{dz}.
\end{eqnarray}
This function is well defined for all complex values of  $z$. Its relation to the original $\Gamma$-function can be obtained by considering $t\to \infty$ limit: 
\begin{eqnarray}
\lim_{t\to \infty}\Gamma_{\mu,t}(z+1)=\mu^{-z}\Gamma(z)
\end{eqnarray}
Also, we note here that this function can be related to the matrix elements of the representations of $G$ (see Section 2) if we would consider Gaussian measure instead of Lebesgue measure on the real line. 
\\

\noindent {\bf 5.3. Laplace transform and representation of the group action via ordinary $\Gamma$-function.} 
In the previous two subsections we defined the analogue of the Fourier transform for the Wiener measure and then related the action of the element of the semigroup $\hat{G}_+$ in the representation $\rho^l_{\lambda,k}$, where $\lambda(u)<0$,  with the $\hat \Gamma$-functional. In this section we relate this action to the ordinary $\Gamma$-function: in order to do that we consider the bilateral Laplace transform with respect to $x_0$ variable. Therefore, we need to restict the representation to the subspace of $\hat{D}\subset H_l$ which consists of functionals  $f(x,x_0)$, such that $f(x,x_0)$ is an infinitely many times differentiable function with respect to $x_0$ with a compact support. Then, if $f(x,z)\in \mathcal{L} \hat D$, the result is as follows:
\begin{eqnarray}
&&\mathcal{L}\rho^l_{\lambda,k}(e^{\alpha},b,s) \mathcal{L}^{-1}f(x,z)=\nonumber\\
&&\frac{1}{2\pi}\int_{\mathbb{R}}e^{ip_0z}
e^{is}e^{-\frac{1}{4t}\int_0^{2\pi}\alpha'(u)\alpha'(u)}e^{-\frac{1}{2t} \int_0^{2\pi}\alpha'(u)\ud x(u)}e^{ik\int_0^{2\pi}\alpha(u)dx(u)}\nonumber\\
&&e^{\int_0^{2\pi} \lambda(u)b(u)e^{x(u)+p_0}\ud u}\int_{\mathbb{R}+i0} e^{-i(p_0+\alpha_0)v} f(x+\t\alpha,v)dvdp_0.
\end{eqnarray}
Since as we assumed $\lambda(u), b(u)>0$, then one can make a change of variables: $p_0\to \eta=p_0+\log(-\int_0^{2\pi} \lambda(u)b(u)e^{x(u)})du$. Then we have:
\begin{eqnarray}\label{act4}
&&\mathcal{L}\rho^l_{\lambda,k}(e^{\alpha},b,s)\mathcal{L}^{-1}f(x,z)=\nonumber\\
&&e^{is}e^{-\frac{1}{4t}\int_0^{2\pi}\alpha'(u)\alpha'(u)}e^{-\frac{1}{2t} \int_0^{2\pi}\alpha'(u)\ud x(u)}e^{ik\int_0^{2\pi}\alpha(u)dx(u)}\nonumber\\
&&\int_{\mathbb{R}+i0} K_{\lambda,\alpha_0,b}(z-v, x(u))f(x+\t \alpha,v)dv,
\end{eqnarray}
 where
\begin{eqnarray}\label{intop}
&& K_{\lambda,\alpha_0,b}(z-v, x)=\nonumber\\
&&  \frac{1}{2\pi}e^{-i\alpha_0w}e^{(iv-iz)\log(-\int_0^{2\pi} \lambda(u)b(u)e^{x(u)})}
\frac{1}{2\pi}\int_{\mathbb{R}} e^{i\eta(z-v)} e^{-e^\eta} d\eta=\nonumber\\
&& \frac{e^{-iz\alpha_0}}{2\pi}(e^{-\alpha_0}\int_0^{2\pi} \lambda(u)b(u)e^{x(u)})^{iv-iz}\Gamma(i(z-v)).
\end{eqnarray}
Similar to the case of $ax+b$ group, one can extend the action of the group elements from $\hat{G}_{+}$ to all complex square integrable $\lambda$ and $b$, such that $\int_0^{2\pi} \lambda(u)b(u)e^{x(u)}du$ is not negative.

Therefore, we have the following theorem.\\

\noindent{\bf Theorem 5.2.} {\it The  action of the element $\mathcal{L}\rho^l_{\lambda,k}(e^{\alpha},b,s)\mathcal{L}^{-1}$ on $\mathcal{L}\hat{\mathcal{D}}$ can be expressed via the integral operator \rf{act4}, with the kernel proportional to the $\Gamma$-function \rf{intop}, as long as ${\rm arg}(\int_0^{2\pi} \lambda(u)b(u)e^{x(u)})<\pi$. }

\section{ Open questions}
\noindent{\bf 6.1. ''Vertex algebra'' related to the representations $\rho^l_{\lambda,k}$. }
The relation between highest weight representations of affine Lie algebras and vertex algebras is well-known (see e.g. \cite{efk}, \cite{benzvi}). Namely, the highest weight modules of the affine algebra are modules over certain vertex operator algebra. 

In Section 3 we described the Lie algebra action in the representation space of $\hat G$. We have two mutually commuting Heisenberg algebras and another object, which reminds one about exponential operators from 
the famous Frenkel-Kac-Segal construction. 
It is clear, that in our case, one has to modify in a certain way the axioms of the vertex algebra in 
order to describe the action of the  relations between generators. 
In such a way this will give us a new object, which will be studied in a separate article. 
\\

\noindent{\bf 6.2. Intertwining operators and possible tensor category for 
$\rho^l_{\lambda,k}$ representations.} 
It is known \cite{ivan} that the unitary representations of $ax+b$ group are closed under tensor product. As we discussed in Section 2 there are three types of ''simple'' objects in the category of unitary representations of $G$. It appears that their tensor products decompose as direct integrals of these ''simple'' objects. 

Using the lessons we learned from the study of the representations of the reductive Lie algebras and their affine analogues, we hope that the same structures would appear in the case of $\hat{G}$ with several modifications. Namely, we expect (this conjecture is due to I.B. Frenkel) to have the braided tensor category, where the braiding is related to the value of the central charge. As well as in the affine case we hope to have a differential equation governing the intertwiners, i.e. the analogue of the Knizhnik-Zamolodchikov equation (see e.g. \cite{efk}).

All this will make possible another intriguing relation.  It is known (this fact is due to Kazhdan and Lusztig) 
that the braided tensor categories of the representations of affine Lie algebras and the quantum algebras, associated with the same simple Lie algebra, are equivalent. 
It is also known that there exists a proper quantum version of the $ax+b$ group, which is called quantum pane. One can expect that the resulting braided tensor category of the representations of this object, studied in \cite{ivan}, \cite{hyunkyu} is related to the possible braided tensor category of the representations of $\hat G$.

\section{Acknowledgements}
I am very grateful to I.B. Frenkel for bringing my attention to this problem and many discussions on the subject. I gained a lot from the remarks of H. Garland and R. Raj. I am grateful to N. Tecu for suggesting several references on the infinite-dimensional analysis during the initial phase of this project. I would like to thank the organizers 
of the Simons Workshop 2010, where this work was partly done.

\end{document}